\documentclass[preprint]{elsarticle}

\usepackage{lineno}
\usepackage{todonotes}
\modulolinenumbers[5]

\journal{Topology and Its Applications}
\usepackage{amscd,amsmath}
\usepackage{amssymb,amsmath,latexsym,color,enumerate,tikz,dsfont}
\usepackage[all]{xy}       
\usepackage{mathabx}
\usepackage{mathrsfs}
\usepackage{stmaryrd}
\usepackage{setspace}
\usepackage[makeroom]{cancel}

\usepackage{hyperref}
\hypersetup{
 unicode=false, 
 pdftoolbar=true, 
 pdfmenubar=true, 
 pdffitwindow=false, 
 pdfstartview={FitH}, 
 pdftitle={comments32}, 
 pdfauthor={JGG,UH,TK}, 
 pdfsubject={points}, 
 pdfcreator={TeXShop}, 
 pdfproducer={Producer}, 
 pdfkeywords={keywords}, 
 pdfnewwindow=true, 
 colorlinks=true, 
 linkcolor=altncolor, 
 citecolor=green, 
 filecolor=blue, 
 urlcolor=cyan 
}

\usepackage{longtable}

\usepackage{floatrow}

\usepackage{amsthm}

\def\frm{{\sf Frm}}

\def\dfrm{{\sf dFrm}}

\def\R{{\mathbb{R}}}

\def\2{\text{\bf 2}}

\def\dFrm{{\mathsf{dFrm}}}

\newcommand{\twoheaddownarrow}{{\rlap{\rlap{$\ $}\raise .25ex\hbox{$\downarrow$}}\raise-.25ex\hbox{$\downarrow$}}}
\newcommand{\twoheaduparrow}{{\rlap{\rlap{$\ $}\raise .25ex\hbox{$\uparrow$}}\raise-.25ex\hbox{$\uparrow$}}}

 \newcommand{\tbigcup}{\mathop{\textstyle \bigcup }}
 \newcommand{\tbigvee
}{\mathop{\textstyle \bigvee }}
 \newcommand{\tbigwedge
}{\mathop{\textstyle \bigwedge }}

\newcommand{\dirvee}{\mathop{\setminus\kern-.85ex\nnearrow}}

\DeclareMathAccent{\widetriangle}{\mathord}{largesymbols}{"E6}

\newcommand{\dva}{\texorpdfstring{$d$}{d}} 
\newcommand{\Dva}{\texorpdfstring{$D$}{d}} 

\newtheorem{theorem}{Theorem}[section]
\newtheorem{proposition}[theorem]{Proposition}
\newtheorem{lemma}[theorem]{Lemma}
\newtheorem{corollary}[theorem]{Corollary}

\theoremstyle{definition}
\newtheorem{definition}[theorem]{Definition}

\newtheorem{examples}[theorem]{Examples}

\newtheorem{notation}[theorem]{Notation}

\theoremstyle{remark}

\newtheorem{remark}[theorem]{Remark}
\newtheorem{remarks}[theorem]{Remarks}

\DeclareMathOperator{\con}{\sf con}
\DeclareMathOperator{\tot}{\sf tot}

\DeclareMathOperator{\dsup}{\mathcal D}

\newcommand{\dirsup}{\mathrel{\setlength{\unitlength}{.7em}\raisebox{-.2em}%
    {\begin{picture}(1,1.5)\put(.5,0){\line(-1,3){.48}}
    \put(.5,0){\vector(1,3){.5}}\end{picture}}}}

\begin{document}

\begin{frontmatter}

\title{On Isbell's Density Theorem\\
for bitopological pointfree spaces I}

\author[CU]{M. Andrew Moshier}
\ead{moshier@chapman.edu}
\author[UPV-EHU]{Imanol~Mozo~Carollo\corref{mycorrespondingauthor}}
\ead{imanol.mozo@ehu.eus}
\author[CU]{Joanne~Walters-Wayland}
\ead{joanne@waylands.com}

\address[CU]{CECAT, Department of Mathematics \& Computer Science, Chapman University,\\ Orange CA 92866, U.S.A.} 
\address[UPV-EHU]{Departamento de Matem\'aticas, Universidad del Pa\'{\i}s Vasco UPV/EHU,\\Apdo. 644, 48080 Bilbao, Spain}

\begin{abstract}
This paper addresses dense sub-objects for point-free bitopology in terms of \dva-frames and provides several examples. We characterize extremal epimorphisms in \dva-frames and show that a smallest dense one always exists, establishing a proper analogue of Isbell's Density Theorem for \dva-frames. Further we explore certain questions about the functoriality of assigning the smallest dense sub-object to each pointfree bitopological space.
\end{abstract}


\end{frontmatter}

\section{Introduction}

Isbell's Density Theorem \cite{I72}, one of the earliest and most telling results of point-free topology --- that every point-free space has a smallest dense subspace --- has nothing to do with the axiom of choice, and everything to do with what ``subspace'' means. 
This shows that point-free topology is not merely an attempt to duplicate point-set topology without the axiom of choice.
In a point-free setting ``sub''  may not mean ``induced on a subset of points.'' Instead, it has to do directly with the behavior of point-free topological embeddings.
Isbell's result, reviewed below, shows clearly that point-free subspaces behave very differently from pointed subspaces.

The reader will be familiar with the basics. 
A point-free topology is a frame --- a complete lattice in which finite meets distribute over all joins. 
A point-free continuous map is the right adjoint to a frame homomorphism --- a function that preserves finite meets and all joins.
So the study of point-free topology generally translates to the study of frames and frame homomorphisms, bearing in mind that a property of frame homomorphisms is dual to the corresponding property of point-free continuous maps.
In this setting, a point-free topological embedding translates to an extremal epimorphism in frames. 
So Isbell is really telling us that there is a smallest dense extremal epimorphism from any frame.

There are two alternative ways to study bitopology in a point-free way, namely biframes and \dva-frames. Each of these gives rise to its own concept of subspace.
The more studied alternative is biframes \cite{BBH83}.
Here a point-free space is thought of in rather close analogy with the pointed variant.
Namely, a biframe consists of a frame $L_0$ and two subframes $L_1$ and $L_2$ so that $L_0$ is generated by $L_1\cup L_2$. 
The idea is that $L_1$ and $L_2$ represent the two topologies of a point-free bitopological space, while $L_0$ represents how they are joined.
Sometimes $L_0$ is called the \emph{ambient} frame of the biframe.
A biframe homomorphism is simply a frame homomorphism of the ambient frames that preserves membership in the two subframes.
So it can be thought of as three frame homomorphisms that commute with the subframe embeddings in the obvious way.


The second alternative, the category of \dva-frames, is the main topic of this paper.
Here a point-free bispace is represented by two frames $(L_-,L_+)$ and 
two relations capturing how they interact.
The first relation, denoted by $\con\subseteq L_+\times L_-$ (consistency), is meant to capture the situation where an ``open'' of $L_+$ and an ``open'' of $L_-$ are ``disjoint''. Of course, in this setting  points are not considered,
so the quote signs are needed.
The second relation, denoted by $\tot\subseteq L_-\times L_+$ (totality), is meant to capture the situation where an ``open'' of $L_-$ and an ``open'' of $L_+$ together ``cover''.

A \dva-frame homomorphism is simply a pair of frame homomorphisms  $(f_-,f_+)$ that together preserve $\con$ and $\tot$.
Density is a fairly obvious adaptation of the familiar notion. 
Namely, a \dva-frame homomorphism is dense if and only if it reflects the $\con$ relation.

Extremal \dva-frame epimorphisms are characterized in this paper, and shown to form a lattice -- but quite a long way from a frame (as happens for  extremal frame epimorphisms).
Indeed, the lattice of extremal epimorphisms does not even have to be modular.

We establish a proper analogue of Isbell's Density Theorem for \dva-frames: every \dva-frame has a smallest dense extremal epimorphism.  In frames, the image of the smallest dense extremal epimorphism is Boolean.
Indeed, ``Booleanization'' is the standard name for this in frames. In \dva-frames, however, the image of the smallest dense epimorphism is not necessarily Boolean.
We interpret this as highlighting that frames are quite special in this respect.
One should not expect, in generalizations of frame theory, that images of smallest dense extremal epimorphisms are necessarily Boolean.

Following the proof of the density theorem, we further investigate the construction of the smallest dense extremal epimorphism as \emph{sub-\dva-locales} --- concretely appearing as special subsets of the underlying \dva-frame structures.
The  result, though explicitly given, is not well-behaved in general. 
First, it is not functorial for all \dva-frame homomorphisms. 
We address this by identifying the morphims for which the construction is functorial. Mirroring the language for frames, we call these morphisms skeletal. 
Second, the most faithful analogue of Isbell's Theorem fails in general because the smallest dense sub-\dva-locale is not comprised only of the double pseudocomplements. 
We characterize the case when the analogy with Isbell fits best: when the double pseudocomplements form the smallest dense sub-\dva-locale. We also investigate the co-reflective subcategory of \dva-frames obtained via the general (ill-behaved) case. 

In comparison to these results, Frith and Shauerte \cite{FS11} show that every biframe has a smallest dense onto biframe homomorphism. But even in biframes, onto is not the same as extremal epimorphism. In a sequel paper, we will consider general extremal epimorphisms in biframes and address Isbell's Density Theorem in that setting.

%
%
%
%


\section{Preliminaries and notation}
%
%
For general notions and results concerning frames we refer to Johnstone \cite{PJ82} or the recent
Picado-Pultr \cite{PP12}. The latter is particularly useful for details about sublocales.

%
%

Isbell's Density Theorem states that each frame has a least dense sublocale. Moreover it has a very precise description: recall that a sublocale $S$ of a frame $L$ is said to be \emph{dense} if $0\in S$. The collection of all elements $a=a^{\ast\ast}$ form a complete Boolean algebra, denoted by $\mathfrak{B}(L)$, the \emph{Booleanization} of $L$, which is the least dense sublocale. The associated quotient is given by the frame homomorphism $b_L\colon L\to\mathfrak{B}(L)$ that maps each element $a$ to its double pseudocomplement $a^{\ast\ast}$.

%
%



In what follows it is necessary to work with pairs of frames, the first of which will be indexed with $-$ and the second indexed with $+$. As a visual aid to type checking, for a pair $(L_-,L_+)$, when the distinction is important, lower case Latin letters ($a$, $b$, $c$, \ldots) will be used for elements in $L_-$ and lower case Greek letters ($\varphi$, $\psi$, \ldots) for elements in $L_+$. A pair of maps ($f_-\colon L_-\to M_-$, $f_+\colon L_+\to M_+$) is often written as $f\colon L\to M$.
%
%
For a relation $\mathcal{R}\subseteq L_-\times L_+$,  $(f_-\times f_+)[\mathcal{R}]\subseteq M_-\times M_+$  will be denoted by $f[\mathcal{R}]$ and, respectively,
for a relation $\mathcal{S}\subseteq L_+\times L_-$,
 $(f_+\times f_-)[\mathcal{S}]\subseteq M_+\times M_-$ by $f[\mathcal{S}]$.

A \emph{\dva-frame} $\mathcal{L}$ is a quadruple $(L_-,L_+,\con,\tot)$ where $L_-$ and $L_+$ are frames, and  \emph{consistency} $\con\subseteq L_+\times L_-$ and \emph{totality} $\tot\subseteq L_-\times L_+$ are relations satisfying the axioms: 
\begin{center}
\renewcommand{\arraystretch}{1.2}
\begin{longtable}{ l l } 
($\con$-$\downarrow$) & $\con$ is a lower set of $L_+\times L_-$, \\
\hline
($\con$-$\vee$) & $\varphi\con a$, $\psi\con b$ implies $\varphi\vee \psi \con a\wedge b$, and $0\con 1$,  \\
($\con$-$\wedge$)& $\varphi\con a$, $\psi\con b$ implies $\varphi\wedge \psi \con a\vee b$, and $1\con 0$,\\
&(i.e. $\con$ is a bounded sublattice of $L_+\times L_-^{\sf op}$)
\\
 \hline
($\tot$-$\uparrow$)& $\tot$ is an upper set in $L_-\times L_+$,\\
\hline
($\tot$-$\wedge$)& $a\tot\varphi$, $b\tot\psi$ implies $a\vee b\tot \varphi\wedge\psi$, and $0\tot 1$,  \\
($\tot$-$\vee$)& $a\tot\varphi$, $b\tot\psi$ implies $a\wedge b\tot \varphi\vee\psi$, and $1\tot 0$,\\
&(i.e. $\tot$ is a bounded sublattice of $L_-^{\sf op}\times L_+$)
\\
\hline
($\con$-$\tot$) & $\con{;}\tot$ is contained in the lattice order on $L_+$, \\
&$\tot^{-1}{;}\con^{-1}$ is contained in the lattice order of $L_-$,\\
\hline
($\con$-$\dirsup$)& $\con$ is closed under directed joins.
\end{longtable}
\end{center}

%
%
%

%
%
%
%

Morphisms $f\colon\mathcal{L}\to\mathcal{M}$ between \dva-frames are pairs $(f_-,f_+)$ of frame homomorphisms preserving the relations $\con$ and $\tot$, meaning that $f[\con_{\mathcal{L}}]\subseteq\con_{\mathcal{M}}$ and $f[\tot_{\mathcal{L}}]\subseteq\tot_{\mathcal{M}}$. The resulting category will be denoted by $\mathsf{dFrm}$. For more details refer to \cite{JM06}.

\begin{remarks}\label{exampledframe}
\begin{enumerate}[(1)]
\item For any two frames $L$ and $M$, define
\[
\varphi\con_0 a\iff \varphi=0\text{ or }a=0\quad\text{and}\quad a\tot_0\varphi\iff a=1\text{ or }\varphi=1.
\]
Then $(L,M,\con_0,\tot_0)$ is a \dva-frame if and only if neither $L$ nor $M$ is the trivial frame, $\boldsymbol{1}$, or both are $\boldsymbol{1}$. Indeed, if, say, $L=\boldsymbol{1}$ then one has $1_{M}\con_0 1\tot_0 0_{M}$, thus $1_{M}=0_{M}$ by ($\con$-$\tot$).
We denote a \dva-frame with minimal $\con$ and $\tot$ as $L.M$.

\item For any frame $L$, define $\con$ and $\tot$ be relations on $L$ given by
\[
\varphi\con a\iff \varphi\wedge a=0\quad\text{and}\quad a\tot\varphi\iff a\vee\varphi=1.
\]
Then $\textit{Sym}(L)=(L,L,\con,\tot)$ is obviously a \dva-frame.
\end{enumerate}
\end{remarks}

\section{Extremal epimorphisms in \texorpdfstring{$\mathsf{dFrm}$}{dFrm}}

%
%

In this section, we characterize extremal epimorphisms and establish a concrete extremal epi-mono factorization system for \dva-frames.
Jakl, Jung and Pultr in \cite{JJP17} establish an extremal epi-mono factorization via a very different method that goes through a co-reflection from a larger category into \dva-frames.
We encourage the reader to compare these methods because we expect them both to  be useful. Further we would like to acknowledge that the  first two results in this section have already been published before with essentially the same proofs. Proposition \ref{monos} is proved in the same way as \cite[Proposition 2.8]{JJP17} and \cite[Proposition 3.3.11]{J18}. Lemma \ref{conScott} is proved in the same way as \cite[Lemma 5]{JJ17} and \cite[Lemma 3.2.1]{J18}. We keep these results in this paper and proofs here to make the paper as self-contained as possible.

To understand extremal epimorphisms, we first need to understand monomorphisms.
Recall that in frames a monomorphism is just a one-one homomorphism. The obvious extension of this notion works for $\dFrm$.

\begin{proposition}\label{monos}
Monomorphisms in $\dfrm$  are precisely those \dva-frame homomorphisms $f$ such that $f_-$ and $f_+$ are monomorphisms in $\frm$.
\end{proposition}
\begin{proof}
If $f_-$ and $f_+$ are one-one then $f$ is obviously a monomorphism: if $f_-(g_-(a))=f_-(h_-(a))$ and $f_+(g_+(\varphi))=f_+(h_+(\varphi))$ then $g_-(a)=h_-(a)$ and $g_+(\varphi)=h_+(\varphi)$. Now, let $f\colon \mathcal{L}\to \mathcal{M}$ be a \dva-frame homomorphism such that, say $f_-$ is not one-one. Then there exists $a,b\in L_-$, $a\neq b$ such that $f_-(a)=f_-(b)$. 
Now consider the \dva-frame $\mathbf{3}.\mathbf{2}$ and the \dva-frame homomorphisms $g,h\colon (\mathbf{3}.\mathbf{2})\to \mathcal{L}$ such that $g_-(c)=x$ and $h_-(c)=y$ (note that $g$ and $h$ obviously exist). Then $fg=fh$ while $g\neq h$.
\end{proof}

The notion of an extremal epimorphism in this category is slightly more complicated; simply taking images works for the frames and the tot relation, but there is no \emph{a priori} reason for the image of the $\con$ relation to satisfy  $(\con$-$\dirsup)$. To address this issue, consider the following:
For a subset $A$ of a frame $(X,\leq)$ let $\dsup(A)$ denote the set of all joins of directed subsets of $A$. If $A$ is a down set, then so is $\dsup(A)$. Indeed, let $x\in L$ such that $x\leq \tbigvee_{i\in I}a_i$ where $\{a_i\}_{i\in I}$ is a directed set of $A$. Then $\{a_i\wedge x\}_{i\in I}$ is a directed set of $A$ and $x=\tbigvee_{i\in I}(a_i\wedge x)$. 

Let $A$ be a down set of $L$ and for ordinals $\alpha$ set
\[
\begin{aligned}
\dsup^0(A)&=A,\quad
\dsup^{\alpha+1}(A)=\dsup(\dsup^\alpha(A)),\quad\text{and}\\
\dsup^\lambda(A)&=\dsup(\tbigcup_{\alpha<\lambda}\dsup^\alpha(A))\quad\text{for limit }\lambda.
\end{aligned}
\]
For a sufficiently large $\alpha$ one has $\dsup^{\alpha+1}(A)=\dsup^\alpha(A)$ and this coincides with $\overline{A}$, the Scott closure of $A$.

\begin{lemma}\label{conScott}
If any relation $\mathcal{R}\subseteq L_+\times L_-$ on frames $L_+$ and $L_-$, satisfies $(\con$-$\vee)$ or  $(\con$-$\wedge)$ then, for any ordinal $\alpha$, so does $\mathcal{D}^\alpha(\mathcal{R})$.
\end{lemma}

\begin{proof} If $\mathcal{R}$ satisfies ($\con$-$\vee$), then given two directed sets $\{(\varphi_i,a_i)\}_{i\in I}$ and $\{(\psi_j,b_j)\}_{j\in J}$ contained in $\mathcal{R}$, the set $\{(\varphi_i\vee\psi_j,a_i\wedge b_j)\}_{i\in I, j\in J}$ is also a directed and contained in $\mathcal{R}$. Since
\[
\tbigvee_{i\in I}\varphi_i\vee\tbigvee_{j\in J}\psi_j=\tbigvee_{i\in I, j\in J}(\varphi_i\vee \psi_j)\quad\text{and}\quad\tbigvee_{i\in I}a_i\wedge\tbigvee_{j\in J}b_j=\tbigvee_{i\in I, j\in J}(a_i\wedge b_j),
\]
$\mathcal{D}(\mathcal{R})$ satisfies ($\con$-$\vee$).
Further, if $\dsup^\alpha(\mathcal{R})$ satisfies ($\con$-$\vee$) for all $\alpha<\lambda$ for some ordinal $\lambda$, then $\tbigcup_{\alpha<\lambda}\dsup^\alpha(\mathcal{R})$ also satifies ($\con$-$\vee$). If $(\varphi,a),(\psi,b)\in \tbigcup_{\alpha<\lambda}\dsup^\alpha(\mathcal{R})$, there exists $\alpha_1,\alpha_2<\lambda$ such that $(\varphi,a)\in\dsup^{\alpha_1}(\mathcal{R})$ and $(\psi,b)\in\dsup^{\alpha_2}(\mathcal{R})$. If, say, $\alpha_1\leq \alpha_2$ one has $(\varphi, a)\in\dsup^{\alpha_2}(\mathcal{R})$ and consequently, $(\varphi\vee\psi, a\wedge b)\in \dsup^{\alpha_2}(\mathcal{R})$. One can check ($\con$-$\wedge$) dually.
\end{proof}

\begin{lemma}
Let $\mathcal{L}$ and $\mathcal{M}$ be \dva-frames and $f\colon \mathcal{L}\to \mathcal{M}$ be a \dva-frame homomorphism.
Then
\[
f(\mathcal{L})=(f_-(L_-),f_+(L_+), \overline{f[\con_{\mathcal{L}}]},f[\tot_{\mathcal{L}}])
\]
is a \dva-frame.
\end{lemma}
\begin{proof}
First note that $f_-(L_-)$ is a subframe of $M_-$, and therefore a frame. Likewise $f_+(L_+)$ is a frame. In order to check ($\tot$-$\vee$), let $a,b\in f_-(L_-)$ and $\varphi, \psi\in f_+(L_+)$ such that
\[
a\tot_{f(L)}\varphi\, \text{and}\, b\tot_{f(L)}\psi.
\]
Then there exist $a^\prime,b^\prime\in L_-$ and $\varphi^\prime, \psi^\prime\in L_+$ with
\[
a= f_-(a^\prime),\, b= f_-(b^\prime),\, \varphi= f_+(\varphi^\prime),\, \text{ and }\, \psi= f_+(\psi^\prime),
\]
such that 
\[
a^\prime\tot_{\mathcal{L}} \varphi^\prime\, \text{ and }\, b^\prime\tot_{\mathcal{L}} \psi^\prime.
\]  Then one has
\[
(a^\prime\vee b^\prime) \tot_{\mathcal{L}} (\varphi^\prime\wedge \psi^\prime)
\]
by ($\tot$-$\wedge$). Consequently
\[
(a\vee b)\tot_{f(L)}(\varphi\wedge\psi)
\]
as 
\[
\varphi\vee\psi=f_+(\varphi^\prime)\vee f_+(\psi^\prime)=f_+(\varphi^\prime\vee\psi^\prime)
\]
and
\[
a\wedge b=f_-(a^\prime)\wedge f_-(b^\prime)=f_-(a^\prime\wedge b^\prime).
\]
One can check ($\tot$-$\wedge$) similarly. The same argument may be used to show that  ($\con$-$\vee$) and ($\con$-$\wedge$) are satisfied by $f[\con_{\mathcal{L}}]$, and hence by $\overline{f[\con_{\mathcal{L}}]}$ (Lemma \ref{conScott}).

The axioms ($\tot$-$\uparrow$) and ($\con$-$\downarrow$) hold since $(f_-\times f_+)\colon L_-\times L_+\to f_-(L_-)\times f_+(L_+)$ is an onto frame homomophism and consequently it maps upper sets into upper sets and lower sets into lower sets. For ($\con$-$\downarrow$), further recall that the Scott closure of a down set is again a down set.
By design, $f(\mathcal{L})$  satisfies ($\con$-$\dirsup$).
In order to check ($\con$-$\tot$), first note that $f[\tot_{\mathcal{L}}]\subseteq\tot_{\mathcal{M}}$ and $f[\con_{\mathcal{L}}]\subseteq \con_{\mathcal{M}}$ since $f$ is a \dva-frame homomorphism.
Since $\con_{\mathcal{M}}$ is Scott-closed, $\con_{f(L)}\subseteq\con_{\mathcal{M}}$.
Consequently $\con_{f(\mathcal{L})};\tot_{f(\mathcal{L})}$ is contained in the lattice order on $f_+(L_+)$ and 
$\tot_{f(\mathcal{L})}^{-1};\con_{f(\mathcal{L})}^{-1}$ is contained in the lattice order on $f_-(L_-)$.
\end{proof}

Consider any \dva-frame morphism $f\colon \mathcal{L}\to \mathcal{M}$.
Let $g\colon \mathcal{L}\to f(\mathcal{L})$ be the pair of frame homomorphisms given by restricting to the codomain of $f$ to the images,
and let $j$ be the pair of frame embeddings from $f(\mathcal{L})$ to $\mathcal{M}$.
The following proposition is immediate from the definitions.

\begin{proposition}\label{factor}
 $g\colon \mathcal{L}\to f(\mathcal{L})$ is a \dva-frame homomorphism,  $j\colon f(\mathcal{L})\to \mathcal{M}$ is a monomorphism in $\dfrm$, and $f=j\cdot g$.
\end{proposition}


This factorization of \dva-frame morphisms leads to the following characterization of extremal epimorphisms. 

\begin{proposition}\label{extremalepis}
Extremal epimorphisms in $\dfrm$ are precisely those \dva-frame homomorphisms $f\colon \mathcal{L}\to \mathcal{M}$ such that $f_-$ and $f_+$ are extremal epimorphisms (i.e., surjective homomorphisms) in $\frm$, $\overline{f[\con_{\mathcal{L}}]}=\con_{\mathcal{M}}$, and $f[\tot_{\mathcal{L}}]=\tot_{\mathcal{M}}$.
\end{proposition}
\begin{proof} 
$(\implies)$: If $f$ does not satisfy any of the conditions, there is the decomposition $f=j\cdot g$ from Proposition \ref{factor} with $j$ a non-isomorphic monomorphism.

($\impliedby$): If $f_-$ and $f_+$ are onto, then $f$ is obviously an epimorphism. Further, let 
\[
\xymatrix{
\mathcal{L}\ar[rr]^{f}\ar[dr]_{g}&&\mathcal{M}\\
&\mathcal{N}\ar[ur]_m&
}
\]
 where $m$ is a monomorphism. Then $m_-$ and $m_+$ are onto and, by Proposition \ref{monos}, they are also one-one. Consequently, $m_-$ and $m_+$ are isomorphisms in \frm. Now 
 \[
 f[\tot_{\mathcal{L}}]=m[g[\tot_{\mathcal{L}}]]\subseteq m[\tot_{\mathcal{N}}]\subseteq \tot_{\mathcal{M}}
 \]
  as $m$ and $g$ are \dva-frame homomorphisms. Therefore, if $f[\tot_{\mathcal{L}}]=\tot_{\mathcal{M}}$, it follows that $m[\tot_{\mathcal{N}}]=\tot_{\mathcal{M}}$. Similarly, 
  \[
 f[\con_{\mathcal{L}}]=m[g[\con_{\mathcal{L}}]]\subseteq m[\con_{\mathcal{N}}]\subseteq \con_{\mathcal{M}}
 \]
 and in consequence $\overline{m[\con_{\mathcal{N}}]}=\con_{\mathcal{M}}$ if $\overline{f[\con_{\mathcal{L}}]}=\con_{\mathcal{M}}$. However, note that $m[\con_{\mathcal{N}}]$ is Scott-closed since $m_+\times m_-$ is a frame isomorphism. Consequently $m$ is an isomorphism.
\end{proof}

Proposition~\ref{extremalepis} shows that Proposition~\ref{factor} yields the anticipated extremal epi-mono factorization system for \dva-frames. 

\section{Sub-\dva-locales}

Since extremal epimorphisms in $\frm$ are represented by unique sublocales, the study of extremal epimorphisms of \dva-frames $\mathcal{L}$ may be reduced to the study of pairs of sublocales of $(L_-,L_+)$.
Given an extremal epimorphism $f\colon\mathcal{L}\to \mathcal{M}$,
there exist sublocales $S_-$ and $S_+$ of $L_-$ and $L_+$, 
and a pair of isomorphisms $g\colon (M_-,M_+)\to (S_-,S_+)$ with $g\cdot f=q$ where $q\colon (L_-,L_+)\to (S_-,S_+)$ is the pair of homomorphisms determined by $(S_-,S_+)$,
that is,
\[
q_{S_-}(a)=\tbigwedge\{s\in S_-\mid s\geq a\}\quad\text{and}\quad q_{S_+}(\varphi)=\tbigwedge\{\psi\in S_+\mid \psi\geq\varphi\}
\]
for each $a\in L_-$ and $\varphi\in L_+$.
The pair of frames $(S_-,S_+)$ can be endowed with $\con_{\mathcal{S}}$ and $\tot_{\mathcal{S}}$ relations induced from those of $\mathcal{M}$ by the isomorphisms $g$ in such a way that $\mathcal{S}=(S_-,S_+,\con_{\mathcal{S}},\tot_{\mathcal{S}})$ is a \dva-frame and $q\colon\mathcal{L}\to \mathcal{S}$ given by $q=(q_{S_-},q_{S_+})$ is an extremal epimorphism. 
Note that, by Proposition \ref{extremalepis},
$\con_{\mathcal{S}}=\overline{q[\con_{\mathcal{L}}]}$ and $\tot_{\mathcal{S}}=q[\tot_{\mathcal{L}}]$.
In consequence, the \dva-frame $\mathcal{S}$ is completely determined by $S_-$ and $S_+$: given a pair of sublocales $S_-$ and $S_+$ there is at most one consistency relation and one totality relation that make $q$ an extremal epimorphism.

\begin{notation}
For a \dva-frame $\mathcal{L}$ and sublocales $S_-$ and $S_+$ of $L_-$ and $L_+$, let
\[
\con_{\mathcal{L}}\cap S=\con_{\mathcal{L}}\cap (S_+\times S_-)\quad\text{and}\quad\tot_{\mathcal{L}}\cap S=\tot_{\mathcal{L}}\cap(S_-\times S_+))
\]
and  call them the \emph{restriction} of $\con_{\mathcal{L}}$ and $\tot_{\mathcal{L}}$ to $S$.
\end{notation}

\begin{lemma}\label{quotienttot}
For any \dva-frame $\mathcal{L}$ and sublocales $S_-$ and $S_+$ of $L_-$ and $L_+$, 
$$q[\tot_{\mathcal{L}}]=\tot_{\mathcal{L}}\cap S$$ where $q=(q_{S_-},q_{S_+})$.
\end{lemma}

\begin{proof}
If $a \tot_{\mathcal{L}} \varphi$ then $q_{S_-}(a)\tot_{\mathcal{L}} q_{S_+}(\varphi)$ by  $\tot$-$\uparrow$ and the fact that $q_{S_-}$ and $q_{S_+}$ are inflationary.
Therefore $q[\tot_{\mathcal{L}}]\subseteq\tot_{\mathcal{L}}\cap S$. Since $q_-$ and $q_+$ are idempotent, we conclude that $\tot_{\mathcal{L}}\cap S\subseteq q[\tot_{\mathcal{L}}]$.
\end{proof}

\begin{definition}
A \emph{sub-\dva-locale} $\mathcal{S}$ of a \dva-frame $\mathcal{L}$ is a \dva-frame where $S_-$ and $S_+$ are sublocales of $L_-$ and $L_+$, and $q=(q_{S_-},q_{S_+})\colon \mathcal{L}\to \mathcal{S}$ is an extremal epimorphism.  The collection of all sub-\dva-locales of $\mathcal{L}$ will be denoted by $\mathsf{dS}(\mathcal{L})$ and is naturally ordered as follows:

\[
\mathcal{S}\leq \mathcal{T}\quad \text{if}
\quad S_-\subseteq T_-\text{ and }S_+\subseteq T_+.
\]
Note that this order relation is the restriction to $\mathsf{dS}(\mathcal{L})$ of the natural preorder relation of extremal epimorphisms in $\dfrm$ with common domain: $f_1\colon \mathcal{L}\to \mathcal{M}_1$ is smaller or equal  to $f_2\colon \mathcal{L}\to \mathcal{M}_2$ if $f_1= g\cdot f_2$ for  some \dva-frame homomorphism $g\colon \mathcal{M}_2\to \mathcal{M}_1$.
\end{definition}

\begin{remark}
Jakl, Jung and Pultr in \cite{JJP17} describe quotients by relations of 
 \dva-frames and show that quotients are extremal epimorphisms. Clearly any sub-\dva-locale $\mathcal{S}$ of a \dva-frame $\mathcal{L}$ is a quotient in the sense of Jakl, Jung and Pultr by a pair of relations $R_-$ and $R_+$ of $L_-$ and $L_+$that yield $S_-$ and $S_+$ as frame quotients, respectively.
\end{remark}

\begin{proposition}
 For any \dva-frame $\mathcal{L}$, $\mathsf{dS}(\mathcal{L})$
forms a complete lattice.  
\end{proposition}
\begin{proof}
$1_{\mathsf{dS}(\mathcal{L})}=\mathcal{L}$ is the top and $0_{\mathsf{dS}(\mathcal{L})}=(\{1\},\{1\},\con,\tot)$  (with the only possible  $\con$ and $\tot$) is the bottom. 
Let $\{\mathcal{S}_i\}_{i\in I}$ be a family of sub-\dva-locales of $\mathcal{L}$ and $(q_{{S_i}_-},q_{{S_i}_+})$ their corresponding extremal epimorphisms. Let  $S_-=\tbigvee_{i\in I}{S_i}_-$ and $S_+=\tbigvee_{i\in I}{S_i}_+$ and $q=(q_{S_-},q_{S_+})$.  Since $q_{S_-}$ and $q_{S_+}$ are onto, $q(\con_{\mathcal{L}})$ satisfies ($\con$-$\downarrow$) and consequently so does $\overline{q(\con_{\mathcal{L}})}$. Clearly, $q(\con_{\mathcal{L}})$ satisfies ($\con$-$\vee$) and ($\con$-$\wedge$), thus, by Lemma \ref{conScott}, so does $\overline{q(\con_{\mathcal{L}})}$. In order to check that 
\[
\mathcal{S}=(S_-,S_+,\overline{q[\con_{\mathcal{L}}]},\tot_{\mathcal{L}}\cap S)
\]
is a \dva-frame it only remains to check that it satisfies ($\con$-$\tot$), as ($\con$-$\dirsup$) is obviously satisfied by design. 
For that purpose, note that 
\[
q_{{S_i}_-}(a)=q_{{S_i}_-}(q_{S_-}(a))\, \text{ and }\, q_{{S_i}_+}(\varphi)=q_{{S_i}_+}(q_{S_+}(\varphi))
\]
for each $a\in L_-, \varphi\in L_+$ and $i\in I$, since $S_i\in S$. Further, as $q_{{S_i}_+}\times q_{{S_i}_-}$ preserves directed sets and arbitrary joins, 
\[
q_{\mathcal{S}_i}\left[\overline{q[\con_{\mathcal{L}}]}\right]\subseteq\overline{q_{\mathcal{S}_i}\left[q[\con_{\mathcal{L}}]\right]}=\overline{q_{\mathcal{S}_i}[\con_{\mathcal{L}}]}=\con_{\mathcal{S}_i}.
\]
Accordingly for any $\varphi, \psi\in  S_+$ and $a\in S_-$ such that \[
\varphi\mathrel{ \overline{q[\con_{\mathcal{L}}]}}a\mathrel{\tot_{\mathcal{L}}} \psi,
\]
it follows that 
\[
q_{\mathcal{S}_i}(\varphi)\mathrel{ \con_{\mathcal{S}_i}} q_{\mathcal{S}_i}(a)\mathrel{\tot_{\mathcal{L}}} q_{\mathcal{S}_i}(\psi).
\]
Therefore, for each $i\in I$, $q_{\mathcal{S}_i}(\varphi)\leq q_{\mathcal{S}_i}(\psi)$, by ($\con$-$\tot$), and hence
\[
\varphi=\tbigwedge_{i\in I} q_{\mathcal{S}_i}(\varphi)\leq \tbigwedge_{i\in I} q_{\mathcal{S}_i}(\psi)=\psi.
\]
Consequently $\con_{\mathcal{S}};\tot_{\mathcal{S}}$ is contained in the lattice order of $S_+$. Dually  it can be checked that $\con_{\mathcal{S}}^{-1} ;\tot_{\mathcal{S}}^{-1}$ is contained in the lattice order of $S_-$.
\end{proof}

Meets of extremal epimorphisms in any category are given by colimits. Thus completeness of $\mathsf{dS}(\mathcal{L})$ follows from the fact that $\mathsf{dFrm}$ is a cocomplete category \cite{JJP17}.
The description of joins of extremal epimorphisms is typically more complicated.
So it is worth noting that our proof explicitly shows how to compute joins rather than meets of sub-\dva-locales and uses it to prove the completeness of $\mathsf{dS}(\mathcal{L})$.

\begin{examples}\label{sdlocadibideak}

\begin{enumerate}[(1)]
\item Let $L_-$ and $_+$ be any non-trivial frames. Then, any pair of non-trivial sublocales $S_-\subseteq L_-$ and $S_+\subseteq L_+$ forms a sub-\dva-locale of $L_-.L_+$. Indeed, let $q_{\mathcal{S}}=(q_{S_-},q_{S_+})$. Note that $q[\con_{L_-.L_+}]$ is the minimal $\con$ relation between $S_-$ and $S_+$, and $q[\tot_{\mathcal{L}}]=\tot_{\mathcal{L}}\cap S$ by Lemma \ref{quotienttot}, and is, therefore, the minimal $\tot$ relation on $S_-$ and $S_+$.
In consequence $\mathcal{S}=(S_-,S_+,q[\con_{\mathcal{L}}],q[\tot_{\mathcal{L}}])$ is a \dva-frame and obviously $q$ is an extremal epimorphisms.

Consider a particular case of the previous example.
The diagram in Figure \ref{subdlocales3} shows the lattice $\mathsf{dS}(\boldsymbol{3}.\boldsymbol{3})$ where $\boldsymbol{3}=\{0<c<1\}$. Each sub-\dva-locale is formed by taking the indicated sublocales with minimal $\con$ and $\tot$.

\begin{figure}\caption{The lattice of sub-\dva-locales of $\boldsymbol{3}.\boldsymbol{3}$}\label{subdlocales3}
\[
\xymatrix{&&\boldsymbol{3}.\boldsymbol{3}\ar@{-}[dll]\ar@{-}[dl]\ar@{-}[dr]\ar@{-}[drr]&&\\
\mathfrak{c}(c).\boldsymbol{3}\ar@{-}[d]\ar@{-}[dr]&\boldsymbol{3}.\mathfrak{c}(c)\ar@{-}[dl]\ar@{-}[drr]&&\boldsymbol{3}.\mathfrak{o}(c)\ar@{-}[dll]\ar@{-}[dr]&\mathfrak{o}(c).\boldsymbol{3}\ar@{-}[dl]\ar@{-}[d]\\
\mathfrak{c}(c).\mathfrak{c}(c)\ar@{-}[drr]&\mathfrak{c}(c).\mathfrak{o}(c)\ar@{-}[dr]&&\mathfrak{o}(c).\mathfrak{c}(c)\ar@{-}[dl]&\mathfrak{o}(c).\mathfrak{o}(c)\ar@{-}[dll]\\
&&\boldsymbol{1}.\boldsymbol{1}&&
}
\]
\end{figure}
This example shows that the lattice of sub-\dva-locales of a \dva-frame may not be distributive, as in this case, for instance, 
\begin{align*}
&\mathfrak{o}(c).\mathfrak{o}(c)\vee\big(\boldsymbol{3}.\mathfrak{c}(c)\wedge\boldsymbol{3}.\mathfrak{o}(c)\big)\\
&\qquad =\mathfrak{o}(c).\mathfrak{o}(c)\vee\boldsymbol{1}.\boldsymbol{1}\\
&\qquad =\mathfrak{o}(c).\mathfrak{o}(c)\\\text{whereas}\\
&\big(\mathfrak{o(c)}.\mathfrak{o}(c)\vee\boldsymbol{3}.\mathfrak{c}(c)\big)\wedge\big(\mathfrak{o(c)}.\mathfrak{o}(c)\vee \boldsymbol{3}.\mathfrak{o}(c)\big)\\
&\qquad =\boldsymbol{3}.\boldsymbol{3}\wedge \boldsymbol{3}.\mathfrak{o}(c)\\
&\qquad =\boldsymbol{3}.\mathfrak{o}(c).
\end{align*}

\item The real line can be considered as a \dva-frame presented as follows \cite{KJM11}:
\[
\mathcal{L}\R=(\mathcal{O}\R_l,\mathcal{O}\R_u,\con_\R,\tot_\R)
\]
where $\mathcal{O}\R_l=\{(-\infty,a)\mid a\in \R\}\cup\{\R,\varnothing\}$ is the frame of open sets of the real line endowed with the lower topology, $\mathcal{O}\R_u=\{(a,+\infty)\mid a\in \R\}\cup\{\R,\varnothing\}$ is the frame of open sets of the real line endowed with the upper topology and $\con_\R$ and $\tot_\R$ are defined in the obvious way. Now let
\begin{align*}
S_-&=\{(-\infty,a)\mid a<0\}\cup\{\varnothing,\R\},\\
S_+&=\{(a,+\infty)\mid a\leq 0\}\cup\{\R\},\\
T_-&=\{(-\infty,a)\mid a\geq 0\}\cup\{\R\},\\
T_+&=\{(a,+\infty)\mid a> 0\}\cup\{\varnothing,\R\},\\
R_-&=\{(-\infty,0),\R\},\quad\text{and}\\
R_+&=\{(0,+\infty),\R\}.
\end{align*}

Note that $S_-,T_-, R_-$ are sublocales of $\mathcal{O}\R_l$ and $S_+,T_+, R_+$ are sublocales of $\mathcal{O}\R_u$.
In particular, $S_-$ and $T_+$ are the open sublocales determined by $(-\infty, 0)$ and $(0,+\infty)$, respectively, and $S_+$ and $T_-$ are the closed sublocales determined by $(-\infty,0)$ and $(0,+\infty)$, respectively.
It is straightforward to check that the pairs $(S_-,S_+), (T_-,T_+)$ and $(R_-,R_+)$ determine sub-\dva-locales of $\mathcal{L}\R$.
Now $S\vee T=\mathcal{L}\R$, since $S_-\vee T_-=\mathcal{O}\R_l$ and $S_+\vee T_+=\mathcal{O}\R_u$, and  
$S\wedge R=T\wedge R=0_{\mathsf{dS}(\mathcal{L}\R)}$, since $S_-\wedge R_-=\{1\}=T_+\wedge R_+$. Therefore, 
\[
(S\vee T)\wedge R=R\neq 0_{\mathsf{dS}(\mathcal{L}\R)}=(S\wedge R)\vee (T\wedge R).
\]
This examples shows that $\mathsf{dS}(\mathcal{L})$ may not be distributive, even in the case of regular \dva-frames ($\mathcal{L}\R$ is regular). Recall that in a \dva-frame $L$ the relation $\con{;}\tot$ of $L_+$ and the relation $\tot^{-1}{;}\con^{-1}$ of $L_-$  are called the \emph{rather below} relations of $L_+$ and $L_-$, respectively. A \dva-frame $\mathcal{L}$ is called \emph{regular} if every element of each component frame is the join of the elements rather bellow it.

\item As we have seen, given a \dva-frame $\mathcal{L}$, not all pairs of sublocales of $L_-$ and $L_+$ form sub-\dva-locales.
The following example shows a little more:
for a sublocale, say, $S_-$ of $L_-$, there may not exist any sublocale $S_+$ of $L_+$ such that the pair $(S_-,S_+)$ forms a sub-\dva-locale of $\mathcal{L}$.
Consider the \dva-frame $\mathcal{L}=(\mathcal{O}(\R),\mathcal{P}(\R), \con_{\mathcal{L}},\tot_{\mathcal{L}})$ where $\mathcal{P}(\R)$ is the powerset of $\R$ ordered by inclusion and $\con_{\mathcal{L}}$ and $\tot_{\mathcal{L}}$  are defined in the obvious way.
There is no sub-\dva-locale $\mathcal{S}$ of $\mathcal{L}$ such that $S_-=\mathfrak{B}(\mathcal{O}(\R))$. 
Assume otherwise. Then for each $x\in\R$ one has
\[
\{x\}\con \R\setminus\{x\}\quad\text{and}\quad q_{S_-}\left(\R\setminus\{x\}\right)=\left(\R\setminus\{x\}\right)^{\ast\ast}=1.
\]
Thus $q_{S_+}(\{x\})\con_{\mathcal{S}} 1$ and consequently,  by ($\con$-$\tot$),  $q_{S_+}(\{x\})=0$ for all $x\in \R$. Since $q_{S_+}(\R)=\tbigvee_{x\in\R}q_{S_+}(\{x\})=0$, $S_+$ is the one element frame $\boldsymbol{1}$, a contradiction, as $S_-=\mathfrak{B}(\mathcal{O}(\R))\neq \boldsymbol{1}$.
\end{enumerate}
\end{examples}

\section{Pseudocomplements in \dva-frames}
In this section we present some elementary properties and results regarding pseudocomplements in \dva-frames that will be necessary to describe dense sub-\dva-locales. For this purpose, let $(L_{-},L_+,\con, \tot)$ be a \dva-frame. 

\begin{definition} For each $a\in L_-$ and $\varphi\in L_+$, let
\[
a^\bullet=\tbigvee\{\varphi\in L_+\mid \varphi\con a\}\quad\text{and}\quad\varphi^\bullet=\tbigvee\{a\in L_-\mid \varphi\con a\}.
\]
These are the \emph{\dva-pseudocomplements} of $a$ and $\varphi$, respectively.
\end{definition}
\begin{remark}\label{remarkspseudo}
The terminology is justified as follows: since $\con$ is closed under directed joins and
\[
\{\varphi\in L_+\mid \varphi\con a\}\quad \text{and}\quad \{a\in L_-\mid \varphi\con a\}
\]
are directed sets by ($\con$-$\vee$) and ($\con$-$\wedge$), respectively, it follows that $a^\bullet$ is the largest element of $L_+$ satisfying $a^\bullet\con a$, and similarly for $\varphi^\bullet$.  
\end{remark}

\begin{lemma}
The maps $(a\mapsto a^\bullet)\colon L_-\to L_+$ and $(\varphi\mapsto\varphi^\bullet)\colon L_+\to L_-$ are antitone.
\end{lemma}

\begin{proof}
This follows easily from ($\con$-$\downarrow$).
\end{proof}

\begin{lemma}\label{eqpseudo}For any $a\in L_-$ and $\varphi\in L_+$, 
\[
\varphi\con a\iff a\leq \varphi^\bullet\iff \varphi\leq a^\bullet.
\]
\end{lemma}

\begin{proof}
There are obvious from Remark~\ref{remarkspseudo}.
\end{proof}

\begin{corollary}\label{Galois} The maps $(a\mapsto a^\bullet)\colon L_-\to L_+$ and $(\varphi\mapsto\varphi^\bullet)\colon L_+\to L_-$ form an antitone Galois connection.
\end{corollary}

The following facts follow by the general theory of Galois connections:

\begin{corollary}\label{propspseudo}For any $a\in L_-$, $\{a_i\}_{i\in I}\subseteq L_-$, $\varphi\in L_+$ and $\{\varphi_j\}_{j\in J}\subseteq L_+$, 
\begin{enumerate}[\rm(1)]
\item $a \leq a^{\bullet\bullet}$,
\item $\varphi\leq \varphi^{\bullet\bullet}$,
\item $a^{\bullet\bullet\bullet}=a^\bullet$,
\item $\varphi^{\bullet\bullet\bullet}=\varphi^{\bullet}$,
\item $(\tbigvee_{i\in I}a_i)^\bullet=\tbigwedge_{i\in I}a_i^\bullet$,
\item $(\tbigvee_{j\in J}\varphi_j)^\bullet=\tbigwedge_{j\in J}\varphi^\bullet$.
\end{enumerate}
\end{corollary}

%
%
%
%

\begin{lemma}\label{condoublepseudo}
For  any $a\in L_-$ and $\varphi\in L_+$,
\[\varphi\con a\iff\varphi^{\bullet\bullet}\con a
\iff\varphi\con a^{\bullet\bullet}.
\]
\end{lemma}

\begin{proof} By Lemma \ref{eqpseudo} and Corollary \ref{propspseudo} (4), 
\[
\varphi^{\bullet\bullet}\con a\iff a\leq \varphi^{\bullet\bullet\bullet}=\varphi^\bullet\iff \varphi\con a.
\] 
 Similarly it can be checked that $ \varphi\con a\iff\varphi\con a^{\bullet\bullet}$.
\end{proof}

\begin{remarks}\label{doubledotex}
\begin{enumerate}[(1)]
\item Note that $L_-^{\bullet\bullet}=\{a^{\bullet\bullet}\mid a\in L_-\}$ and $L_+^{\bullet\bullet}=\{\varphi^{\bullet\bullet}\mid\varphi\in L_+\}$ are precisely the subsets of idempotent elements of a Galois connection. Therefore $L_-^{\bullet\bullet}$  and $L_+^{\bullet\bullet}$  are complete lattices in which meets are the same as in $L_-$ and $L_+$, respectively. Also, joins in $L_-^{\bullet\bullet}$ are given by 
$\left(\tbigvee_{i\in I}a_i\right)^{\bullet\bullet}$, and joins in $L_+^{\bullet\bullet}$ are given by $\left(\tbigvee_{j\in J}\varphi_j\right)^{\bullet\bullet}$. 

\item However $L_-^{\bullet\bullet}$ and $L_+^{\bullet\bullet}$ are not frames in general. For instance, let $L_-=\mathcal{O}\R$ be the Euclidean topology on the real line, $L_+$ the cocompact topology on the real line and $\mathcal{L}=(L_-,L_+,\con,\tot)$ be a \dva-frame where $\con$ and $\tot$ are defined in the obvious way. Note that for any unbounded open subset $U$ in $L_-$, one has $U^\bullet=\varnothing$. Then one has
\[
\begin{aligned}
(0,1)\wedge\tbigvee^{L_-^{\bullet\bullet}}_{n\geq 1}(n,n+1)&=(0,1)\cap\left(\tbigcup_{n\geq 1}(n,n+1)\right)^{\bullet\bullet}\\
&=(0,1)\wedge\R=(0,1)
\end{aligned}
\]
while
\[
\tbigvee^{L_-^{\bullet\bullet}}_{n\geq 1}((0,1)\wedge(n,n+1))=\varnothing.
\]
Therefore $L_-^{\bullet\bullet}$ is not a frame.
\end{enumerate}
\end{remarks}

%
%


%
%

\section{The \dva-localic version of Isbell's Density Theorem}

Isbell's Density Theorem has two important features. First, every frame has a least dense sublocale.
This is easy to see because (i) sublocales are closed under arbitrary intersection, 
(ii) the given frame is dense in itself and (iii) intersections of dense sublocales are dense. 
The more surprising, second, aspect is that the least dense sublocale consists exactly of all (double) pseudo-complements, and is therefore Boolean. 
Indeed it is this that leads authors to call it the Booleanization of the frame.
Banaschewski and Pultr \cite{BP96} characterize the frame homomorphisms for which Booleanization is functorial, thus demonstrating Booleanization as the co-reflection of the category frames with skeletal homomorphims into the category of Boolean frames with all frame homomorphisms.

In \dva-frames, even the first feature is not simple. The notion of density  for \dva-frames is stronger than mere density of the two component sublocales. So one can not simply apply Isbell's theorem to the two components of a \dva-frame and expect anything useful to come of it. Instead, a version of Isbell's theorem for \dva-frames requires a necessary and sufficient condition for density of sub-\dva-locales that clearly is closed under arbitrary intersection of sublocales in the two component frames. Then the three point outline of Isbell's proof adapts.

The second feature of Isbell's result is also more subtle for \dva-frames. The dense sub-\dva-locales are those that include all double-pseudocomplements,
but as we have seen in Remark \ref{doubledotex}, the double-pseudocomplements do not always form sublocales.
Even in case they do, they do not form Boolean sublocales.

Recall that a frame homomorphism is \emph{dense} if it reflects $0$:  $f(a)=0 \implies a = 0$. Adapted to \dva-frames, density means instead that $\con$ is reflected.

\begin{definition}
A \dva-frame homomorphism $f\colon \mathcal{L}\to \mathcal{M}$ is said to be \emph{dense} if
\[
f_+(\varphi)\con_{\mathcal{M}} f_-(a)\implies \varphi\con_{\mathcal{L}} a.
\]
Analogously, a sub-\dva-locale $\mathcal{S}$ of a \dva-frame $\mathcal{L}$ is \emph{dense} if the extremal epimorphism $q\colon \mathcal{L}\to \mathcal{S}$ determined by $\mathcal{S}$ is dense, that is, if $\con_\mathcal{S}$  is the restriction of $\con_{\mathcal{L}}$  to $\mathcal{S}$.
\end{definition}

\begin{remark}\label{exampledense} If $f\colon \mathcal{L}\to \mathcal{M}$ is a dense \dva-frame homomorphism, then $f_-$ and $f_+$ are dense frame homomorphisms. Indeed, if, say, $a\in L_-$ is such that $f_-(a)=0$ then $f_+(1)=1\con_{\mathcal{M}} f_-(a)$. If $f$ is dense this implies $1\con_{\mathcal{L}} a$ and consequently $0\tot_{\mathcal{L}} 1\con_{\mathcal{L}} a$. By ($\con$-$\tot$), it follows that $a=0$. That $f_+$ is dense follows dually.
%
%

The converse does not hold. For instance, let $L_-$ be the topology on the set $X=\{a,b,c\}$ with $\{a,b\}$ and $\{a\}$ the only non-trivial open subsets and $L_+$ be the topology on $X$ with $\{b,c\}$ the only non-trivial open subset. Further let $M_-$ be the topology on $X$ with $\{a\}$ the only non-trivial open subset. Let $\mathcal{L}=(L_-,L_+,\con_{\mathcal{L}},\tot_{\mathcal{L}})$ and $\mathcal{M}=(M_-,L_+,\con_{\mathcal{M}},\tot_{\mathcal{M}})$ be \dva-frames where the $\con$ and $\tot$ relations are defined on the obvious way from the set theoretical relations. Let $f_-\colon L_-\to M_-$ be given by $f_-(\{a,b\})=\{a\}=f_-(\{a\}),f_-(\varnothing)=\varnothing$ and $f_-(X)=X$. $f_-$ is obviously a dense homomorphism. Let $f_+$ be the identity map for $L_+$. Then $f=(f_-,f_+)$ is \dva-frame homomorphism given by a pair of dense frame homomorphisms, but $f$ is not a dense homomorphism, as $f_+(\{b,c\})\con_{\mathcal{M}} f_-(\{a,b\})$ while $\{b,c\}\cancel{\con_{\mathcal{L}}}\{a,b\}$.

\end{remark}


\begin{lemma}\label{densepseudo1}
Let $\mathcal{S}$ be a dense sub-\dva-locale of a \dva-frame $\mathcal{L}$ and $q\colon \mathcal{L}\to \mathcal{S}$ the extremal epimorphism determined by $\mathcal{S}$. Then 
\[
q_+(a^\bullet)=a^\bullet\quad\text{and}\quad q_-(\varphi^\bullet)=\varphi^\bullet.
\]
\end{lemma}

\begin{proof}
Since $q$ is inflationary, $a^\bullet\leq q_+(a^\bullet)$.  On the other hand, $q_+(a^\bullet)\con_\mathcal{S} q_-(a)$ since $a^\bullet\con_{\mathcal{L}} a$. But $\mathcal{S}$ is dense, so $q_+(a^\bullet)\con_{\mathcal{L}} q_-(a)$ and, consequently, $q_+(a^\bullet)\con_{\mathcal{L}} a$ by ($\con$-$\downarrow$) and the fact that $a\leq q_-(a)$. By Lemma \ref{eqpseudo}, it follows that $q_+(a^\bullet)\leq a^\bullet$ .
\end{proof}
%
%
\begin{proposition}\label{densepseudo2}
Let $\mathcal{L}$ be a \dva-frame, $S_-$ be a sublocale of $L_-$, and $S_+$ be a sublocale of $L_+$ such that $L_-^{\bullet\bullet}\subseteq S_-$ and $L_+^{\bullet\bullet}\subseteq S_+$. 
Then $(S_-,S_+,\con_\mathcal{S},\tot_\mathcal{S})$ forms a dense sub-\dva-locale of $\mathcal{L}$, where $\con_\mathcal{S}$ and $\tot_\mathcal{S}$ are the restrictions of $\con_{\mathcal{L}}$ and $\tot_{\mathcal{L}}$ to $S_+\times S_-$ and $S_-\times S_+$, respectively.
\end{proposition}

\begin{proof} It is easy to verify that ($\con$-$\downarrow$), ($\tot$-$\uparrow$) and ($\con$-$\tot$) hold for the relations $\con_\mathcal{S}=\con_{\mathcal{L}}\cap S$ and $\tot_\mathcal{S}=\tot_{\mathcal{L}}\cap S$. In order to check ($\tot$-$\wedge$), let $a,b\in S_-$ an $\varphi,\psi\in S_+$ such that $a\tot_\mathcal{S} \varphi$ and $b\tot_\mathcal{S}\psi$. Then 
\[
a\vee_S b\geq a\vee_L b\tot_{\mathcal{L}} \varphi\wedge_L\psi=\varphi\wedge_S\psi.
\]
It follows that $a\vee_S b\tot_{\mathcal{S}} \varphi\wedge_S\psi$ by ($\tot$-$\uparrow$). Dually ($\tot$-$\vee$).

To show ($\con$-$\vee$), let $a,b\in S_-$ and $\varphi,\psi\in S_+$ such that $\varphi\con_\mathcal{S} a$ and $\psi\con_\mathcal{S} b$. Then, by ($\con$-$\vee$), one has $\varphi\vee_L \psi\con_{\mathcal{L}} a\wedge_L b$. By Lemma \ref{condoublepseudo}, $(\varphi\vee_L\psi)^{\bullet\bullet}\con_{\mathcal{L}} a\wedge_L b$. Since $L_-^{\bullet\bullet}\subseteq S_-$, $\varphi\vee_S\psi\leq (\varphi\vee_L\psi)^{\bullet\bullet}$ and thus
\[
\varphi\vee_S\psi\con_{\mathcal{S}} a\wedge_S b=a\wedge_L b.
\]
Analogously for ($\con$-$\wedge$).

In order to check that $\con_{\mathcal{S}}$ is closed under directed joins, let $\{(\varphi_i,a_i)\}_{i\in I}\subseteq\con_{\mathcal{S}}$ be a directed set. Obviously, $\{(\varphi_i,a_i)\}_{i\in I}\subseteq\con_{\mathcal{L}}$, and in consequence
\[
\tbigvee^L_{i\in I}\varphi_i\con_{\mathcal{L}} \tbigvee^L_{i\in I}a_i.
\]
By Lemma \ref{condoublepseudo},
\[
\left(\tbigvee_{i\in I}^L\varphi_i\right)^{\bullet\bullet}\con_{\mathcal{L}} \left(\tbigvee_{i\in I}^L 
a_i\right)^{\bullet\bullet}.
\]
Now, since $L_-^{\bullet\bullet}\subseteq S_-$ and $L_+^{\bullet\bullet}\subseteq S_+$, 
\[
\tbigvee_{i\in I}^S\varphi_i\leq \left(\tbigvee_{i\in I}^L\varphi_i\right)^{\bullet\bullet}\quad\text{and}\quad\tbigvee_{i\in I}^S a_i\leq \left(\tbigvee_{i\in I}^L a_i\right)^{\bullet\bullet}
\]
 and ($\con$-$\downarrow$) gives  that 
 \[
\tbigvee_{i\in I}^S\varphi_i\con_{\mathcal{S}} \tbigvee_{i\in I}^S
a_i.
\]
Therefore $(S_-,S_+,\con_{\mathcal{S}},\tot_{\mathcal{S}})$ forms a sub-\dva-locale of $\mathcal{L}$. Finally, note that it is obviously dense.
\end{proof}

This gives the following characterization of dense sub-\dva-locales: 

\begin{corollary}
Let $\mathcal{S}$ be a sub-\dva-locale of a \dva-frame $\mathcal{L}$. Then $\mathcal{S}$ is dense if and only if $L_-^{\bullet\bullet}\subseteq S_-$ and $L_+^{\bullet\bullet}\subseteq S_+$.
\end{corollary}

\begin{proof}
Necessity follows easily from Lemma \ref{densepseudo1} and sufficiency from Proposition \ref{densepseudo2}.
\end{proof}

The preceding lemmas come together to prove the \dva-frame analogue of Isbell's Density Theorem.
%
%

%
%
%
\begin{theorem}\label{isbell}
Each \dva-frame has a smallest dense sub-\dva-locale.
\end{theorem}

\begin{proof}
For a \dva-frame $\mathcal{L}$, take $\hat{L}_-$ and $\hat{L}_+$ to be the smallest sublocales of $L_-$ and $L_+$ containing $L_-^{\bullet\bullet}$ and $L_+^{\bullet\bullet}$, respectively.
By Proposition~\ref{densepseudo2}, these sublocales induce a dense sub-\dva-locale of $\mathcal{L}$. 
By construction and the corollary above, any other dense sub-\dva-locale must include this.
\end{proof}

We let $\hat{\mathcal{L}}$ denote the smallest dense sub-\dva-locale of $\mathcal{L}$.

\begin{examples}
\begin{enumerate}[(1)]
\item The smallest dense sub-\dva-locale of the \dva-frame of the reals $\mathcal{L}\R=(\mathcal{O}\R_l,\mathcal{O}\R_u,\con_\R,\tot_\R)$ is itself. The smallest dense sub-\dva-locale of $\boldsymbol{3}.\boldsymbol{3}$ in Examples \ref{sdlocadibideak} (1) is also itself.

\item Given any frame $L$ , the \dva-frame $Sym(L)$ has $L_-^{\bullet\bullet}=\hat{L}_-=\mathfrak{B}L$ and $L_+^{\bullet\bullet}=\hat{L}_+=\mathfrak{B}L$.

\item
Recall the \dva-frame $\mathcal{L}$ in Remarks \ref{doubledotex} (2) with $L_-$ the Euclidean topology on the real line, $L_+$ the  cocompact topology on the real line  and $\con$ and $\tot$ given in the obvious way. Note that $L_+\subseteq L_-$. One can easily show that $L_+^{\bullet\bullet}=\mathfrak{B}(L_-)\cap L_+$. Since the assignment $a\mapsto a^{\bullet\bullet}$ determines a nucleus on $L_-$ and its restriction to $L_+$ takes images in $L_+$, this restriction is indeed a nucleus in $L_+$. In consequence $L_+^{\bullet\bullet}$ is a sublocale of $L_+$ and $L_+^{\bullet\bullet}=\hat{L}_+$. Further, one can easily check that
\[
L_-^{\bullet\bullet}=\{a\in\mathfrak{B}(L_-)\mid a\text{ is bounded}\}.
\]
Therefore, as $0\in L_-^{\bullet\bullet}\subseteq\mathfrak{B}(L_-)$, then $\hat{L}_-=\mathfrak{B}(L_-)$.
\end{enumerate}
\end{examples}

\section{Good and bad behavior}

%
Although Theorem~\ref{isbell} is an analogue of Isbell's theorem, the analogy is imperfect for interesting reasons.
In particular, Remark~\ref{doubledotex} shows that $L_-^{\bullet\bullet}$ and $L_+^{\bullet\bullet}$ 
are not always sublocales. 
So $\hat{L}_-$ may properly contain $L_-^{\bullet\bullet}$,
and similarly for $\hat{L}_+$. 
Immediately, this raises the question of 
what conditions ensure that  $L_-^{\bullet\bullet}$ and $L_+^{\bullet\bullet}$ are sublocales.
This also raises the quesiton of what are the conditions in which this ``double pseudocomplement'' construction yields a co-reflection.
To address this, we must also understand which class of \dva-frame morphisms make $\mathcal{L}\mapsto \hat{\mathcal{L}}$ functorial.

%
%

For a \dva-frame $\mathcal{L}$, define a relation  $\sqsubseteq_-$ on $L_-$ by $a\sqsubseteq_- b$ if and only if for all $c\in L_-$, $\varphi\in L_+$, $\varphi\con b\wedge c$ implies $\varphi\con a\wedge c$. 
Evidently, $a\leq b$ implies $a\sqsubseteq_- b$, and $\sqsubseteq_-$ is transitive. 

\begin{lemma}\label{framerel}
Let $\mathcal{L}$ be a \dva-frame. Then
\begin{enumerate}[\rm(1)]
\item for any $\{a_i\}_{i\in I}\subseteq L_-$ and $b\in L_-$, if $a_i\sqsubseteq b$ for all $i\in I$, then $\tbigvee_{i\in I} a_i\sqsubseteq_- b$;
\item for any $a,b_0,b_1\in L_-$, if $a\sqsubseteq_- b_0$ and $a\sqsubseteq_- b_1$, then $a\sqsubseteq_- b_0\wedge b_1$.
\end{enumerate}
\end{lemma}
\begin{proof}
Suppose $a_i\sqsubseteq_- b$ holds for each $i\in I$. Then for any $c$ and $\varphi$ satisfying $\varphi\con b\wedge c$, it is the case that $\varphi\con a_i\wedge c$. Using the axiom $\con$-$\vee$ for any finite subset $J\subseteq I$, $\varphi \con \tbigvee_{j\in J}a_i \wedge c$. So by
$\con$-$\dirsup$, $\varphi\con \tbigvee_{i\in I}a_i\wedge c$.

Suppose that for all $c$ and $\varphi$, (i) $\varphi\con b_0\wedge c$ implies $\varphi\con a\wedge c$,
and (ii) $\varphi\con b_1\wedge c$ implies $\varphi\con a\wedge c$. 
Consider some $c$ and $\varphi$ for which $\varphi\con b_0\wedge (b_1\wedge c)$. 
Then by (i) $\varphi\con a\wedge (b_1\wedge c)$. So by (ii) $\varphi\con a\wedge (a\wedge c)$. 
\end{proof}

\begin{lemma}\label{lemma7.2}
Let $\mathcal{L}$ be a \dva-frame. Then for $a\in L_-$, the following are equivalent:
\begin{enumerate}[\rm(1)]
\item $a\in\hat{L}_-$,
\item $b\sqsubseteq_- a$ implies $b\leq a$ for all $b\in L_-$,
\item $a = \tbigvee\{b\in L_-\ \mid\ b\sqsubseteq_- a\}$.
\end{enumerate}
\end{lemma}
\begin{proof}
Condition (2) obviously implies that $\tbigvee\{b\in L_-\ \mid\ b\sqsubseteq_- a\}\leq a$.
But for any $a\in L_-$, it is the case that $a\sqsubseteq_- a$. So $a\leq \tbigvee\{b\in L_-\ \mid\ b\sqsubseteq_- a\}$. And (3) implies (2), trivially.

The map $\nu_-\colon L_-\to L_-$ given by $a\mapsto \tbigvee\{b\in L_-\ \mid\ b\sqsubseteq_- a\}$ is clearly monotonic.
It is inflationary because $a\sqsubseteq_- a$ holds for all $a\in L_-$. 
Idempotency follows from Lemma \ref{framerel}.
And because $b\sqsubseteq a_0\wedge a_1$ if and only if $b\sqsubseteq_- a_0$ and $b\sqsubseteq_- a_1$, the map $\nu_-$ is a nucleus.
So the elements satisfying (3) form a sublocale $\nu_-(L_-)$. 

Consider $\varphi\in L_+$. We claim that $\varphi^\bullet\in \nu_-(L_-)$. 
Suppose for some $b\in L_-$, it is the case that for all $c\in L_-$, $\psi\in L_+$,
if $\psi\con \varphi^\bullet \wedge c$, then $\psi\con b\wedge c$.
In particular, $\varphi\con \varphi^\bullet\wedge 1$. So $\varphi\con b\wedge 1$. 
Hence $b\leq \varphi^\bullet$. 
This shows that $L^{\bullet\bullet}_-\subseteq \nu_-(L_-)$. 
Hence (1) implies (3).

Suppose $S$ is a sublocale of $L_-$ containing $L_-^{\bullet\bullet}$.
Let $a$ be an element satisfying (2). That is, $b\sqsubseteq_- a$ implies $b\leq a$.
So if for all $c$ and $\psi$, $\psi\con a\wedge c$ implies $\psi\con b\wedge c$, then $b\leq a$.

For each $c$ and $\varphi$, the element $c\to \varphi^\bullet$ belongs to $S$ because $\varphi^\bullet\in L^{\bullet\bullet}_-$ and $S$ is a sublocale. 
Notice that $\varphi\con x\wedge c$ if and only if $x\leq c\to \varphi^\bullet$.
So $a$ is the meet of elements $c\to\varphi^\bullet$ such that $\varphi\con a\wedge c$.
Thus $a\in S$. 
\end{proof}

As a corollary of the previous proof we have the following result.
\begin{lemma}\label{quotientsq}
Let $\mathcal{L}$ be a \dva-frame  and the maps  $\nu_-\colon L_-\to L_-$ given by \[a\mapsto \tbigvee\{b\in L_-\ \mid\ b\sqsubseteq_- a\}.\]
Then \[\hat{L}_-=\nu_-(L_-)\cong L_-/\sqsubseteq_-.\]
\end{lemma}

\begin{lemma}\label{doubledotnuclear}
In a \dva-frame $L$, the following are equivalent:
\begin{enumerate}[\rm(1)]
\item $L_-^{\bullet\bullet}$ is a sublocale of $L_-$;
\item $L_-^{\bullet\bullet} = \hat{L}_-$;
\item $b\leq a^{\bullet\bullet}$ if and only if $b\sqsubseteq_- a$ for all $a,b\in L_-$;
\item $(a\wedge b)^{\bullet\bullet} = a^{\bullet\bullet}\wedge b^{\bullet\bullet}$ for all $a,b\in L_-$;
\item $(a\wedge b)^\bullet = (a^{\bullet\bullet} \wedge b)^\bullet$ for all $a,b\in L_-$. 
\item $ \varphi\con a\wedge b$  implies $\varphi\con a^{\bullet\bullet}\wedge b$ for all $a,b\in L_-$ and all $\varphi\in L_+$;
\item $a\sqsubseteq_- b$ implies $a^{\bullet\bullet}\sqsubseteq_- b$.
\end{enumerate}
\end{lemma}
\begin{proof}
If $L_-^{\bullet\bullet}$ is a sublocale, it is, trivially, the smallest one containing $L_-^{\bullet\bullet}$.
If $L_-^{\bullet\bullet} = \hat{L}_-$, then $a\mapsto a^{\bullet\bullet}$ must be the associated nucleus.
If $a\mapsto a^{\bullet\bullet}$ is a nucleus, then it preserves finite meets.

Suppose $(a\wedge b)^{\bullet\bullet} = a^{\bullet\bullet}\wedge b^{\bullet\bullet}$ holds everywhere. 
Since $a\mapsto a^{\bullet\bullet}$ is monotonic, inflationary and idempotent, 
it is a nucleus, so its image, $L_-^{\bullet\bullet}$, is a sublocale.

Suppose $(a\wedge b)^{\bullet\bullet} = a^{\bullet\bullet} \wedge b^{\bullet\bullet}$ holds everywhere. Then 
\begin{align*}
(a\wedge b)^\bullet &= (a\wedge b)^{\bullet\bullet\bullet}\\
  &= (a^{\bullet\bullet}\wedge b^{\bullet\bullet})^\bullet\\
  &= (a^{\bullet\bullet\bullet\bullet}\wedge b^{\bullet\bullet})^\bullet\\
  &= (a^{\bullet\bullet}\wedge b)^{\bullet\bullet\bullet}\\
  &= (a^{\bullet\bullet}\wedge b)^\bullet
\end{align*}

Suppose $(a\wedge b)^\bullet = (a^{\bullet\bullet}\wedge b)^\bullet$ everywhere. Then
\begin{align*}
(a\wedge b)^{\bullet\bullet} &\leq a^{\bullet\bullet}\wedge b^{\bullet\bullet}\\
&\leq (a^{\bullet\bullet}\wedge b^{\bullet\bullet})^{\bullet\bullet}\\
 &= (a\wedge b)^{\bullet\bullet}
\end{align*}
Evidently, $(a^{\bullet\bullet}\wedge b)^\bullet  \leq (a\wedge b)^\bullet$ in any \dva-frame. The opposite comparison is clearly  equivalent to condition (6).

Condition (7) is simply a direct restatement of (6). Simply note that if (7) holds, then $a^{\bullet\bullet}\sqsubseteq_- a$ follows from $a\sqsubseteq_- a$.
\end{proof}

The relation $\sqsubseteq_+$ defined similarly on $L_+$ satisfies lemmas exactly analogous to Lemmas~\ref{framerel} through \ref{doubledotnuclear}.

Say that a \dva-frame having the properties stated in Lemma~\ref{doubledotnuclear} in both $L_-$ and $L_+$ is \emph{corrigible}. Remark \ref{doubledotex} (2) provides an incorrigible \dva-frame. 

One can think of corrigibility as a localness property.
That is, in any \dva-frame, $a^\bullet=a^{\bullet\bullet\bullet}$ can be rewritten $(a\wedge 1)^\bullet = (a^{\bullet\bullet}\wedge 1)^\bullet$.
So (6) asserts that this holds for any $b$ in place of $1$.

The next question we address is how to characterize the \dva-frame morphisms for which $L\mapsto \hat{L}$ extends to a functor.
It is clear that $\mathcal{L}\mapsto\hat{\mathcal{L}}$ cannot be extended to a functor for arbitrary \dva-frame morphisms, because if it could be, then by restricting to symmetric \dva-frames, Isbell's theorem would yield a co-reflection of Boolean frames in the category of frames.

For \dva-frame morphism $f\colon \mathcal{L}\to \mathcal{M}$,  say that $f$ is \emph{skeletal} if  $f_-$ and $f_+$ preserve the relations $\sqsubseteq_-$ and $\sqsubseteq_+$. 
Note that if $\mathcal{M}$ is corrigible, this definition takes the more recognizable form of $f_-(a^{\bullet\bullet}) \leq f_-(a)^{\bullet\bullet}$ and $f_+(\varphi^{\bullet\bullet}) \leq f_+(\varphi)^{\bullet\bullet}$ \cite{BP96}. 

\begin{lemma}\label{hatfunctoriality}
For \dva-frame morphism $f\colon\mathcal{L}\to\mathcal{M}$, define $\hat{f}\colon\mathcal{L}\to\mathcal{M}$ by 
\begin{align*}
\hat{f}_-(a) &=\tbigvee\{b\in M_-\ \mid\ b\sqsubseteq_- f_-(a)\}\\
\hat{f}_+(\varphi) &= \tbigvee\{\psi\in M_+\ \mid\ \psi\sqsubseteq_+ f_+(\varphi)\}
\end{align*}
Then $\hat{f}$ co-restricts to $\hat{\mathcal{M}}$. Moreover, for skeletal morphisms, $f\mapsto \hat{f}$ is functorial. 
\end{lemma}
\begin{proof}
The formula defining $\hat{f}_-$ clearly sends 
any element of $L_-$ to an element of $\hat{M}_-$.
Likewise for $\hat{f}_+$. 
So by restriction $\hat{f}_-$ is a function from $\hat{L}_-$ to $\hat{M}_-$, and similarly for $\hat{f}_+$.

The identity on a \dva-frame $\mathcal{L}$ is obviously skeletal, and its restriction to $\hat{\mathcal{L}}$ is identity.

Let $f\colon \mathcal{L}\to \mathcal{M}$ be a \dva-frame morphism and $g\colon \mathcal{M}\to \mathcal{N}$ a skeletal \dva-frame morphism.
Define $\nu_-^{\mathcal{N}}(b) = \tbigvee\{c\in N_-\ \mid\ c\sqsubseteq_- b\}$ and likewise $\nu_-^{\mathcal{M}}(b) = \tbigvee\{c\in M_-\ \mid\ c\sqsubseteq_- b\}$. 
As we have checked in Lemma \ref{lemma7.2}, $\nu_-^{\mathcal{N}}$ is monotonic, idempotent and inflationary. And
\begin{align}
\hat{g}_-(\hat{f}_-(a)) &= \nu_-^N(g_-(\nu_-^M(f_-(a)))\\
\widehat{g\circ f}_-(a) &= \nu_-^N(g_-(f_-(a)))
\end{align}
Evidently $\widehat{g\circ f}_-(a)\leq \hat{g}_-(\hat{f}_-(a))$ because $\nu^M$ is inflationary and both $\nu^N$ are $g_-$ are monotonic. Conversely, 
$\widehat{g\circ f}_-(a)\geq \hat{g}_-(\hat{f}_-(a))$ because $g$ is skeletal and $\nu_-^{\mathcal{N}}$ is idempotent.
\end{proof}

Say that a \dva-frame $L$ is a \emph{double negation \dva-frame} if $a^{\bullet\bullet}=a$ and $\varphi^{\bullet\bullet}=\varphi$ for all $a$ and $\varphi$. 

\begin{remark}
Heyting algebras (hence frames, in particular) have the property that the law of double negation and the law of the excluded middle are equivalent, both characterizing Boolean algebras among Heyting algebras. 
In \dva-frames the situation is quite different.
Call a \dva-frame an \emph{excluded middle \dva-frame} if $a\tot a^\bullet$ and $\varphi^\bullet\tot \varphi$ for all $a$ and $\varphi$.
Then an excluded middle \dva-frame is a double negation \dva-frame. 
That is, $a^\bullet\con a^{\bullet\bullet}$ always holds. So $a \tot a^\bullet$ implies $a^{\bullet\bullet}\leq a$.
On the other hand, double negation does not imply excluded middle. 
The \dva-frame of the reals is a counter example.
\end{remark}

\begin{lemma}\label{doubleneglemma}
A double negation \dva-frame is corrigible. Moreover, in such a frame $a \sqsubseteq_- b$ implies $a\leq b$
and $\varphi\sqsubseteq_+ \psi$ implies $\varphi\leq \psi$. Consequently, all \dva-frame morphisms from double negation \dva-frames are skeletal.
\end{lemma}
\begin{proof}
If $a^{\bullet\bullet} = a$, then trivially $a\sqsubseteq_- b$ implies $a^{\bullet\bullet}\sqsubseteq_- b$. 
Suppose $a\sqsubseteq_- b$. Then \[\forall d\in L_-\forall\varphi\in L_+, \varphi\con b\wedge d\Rightarrow \varphi\con a\wedge d.\]
In particular, setting $d=1$, $\varphi\con b$ implies $\varphi\con a$. So $b^\bullet\leq a^\bullet$. 
Hence in a double negation frame, $a=a^{\bullet\bullet}\leq b^{\bullet\bullet} = b$.
\end{proof}

\begin{theorem}\label{doublenegtionCoreflection}
  The category of double negation \dva-frames is co-reflective in the category of corrigible \dva-frames with skeletal morphisms.
\end{theorem}
\begin{proof}
From Lemma~\ref{doubleneglemma}, $\hat{\hat{\mathcal{L}}}$ is isomorphic to $\hat{\mathcal{L}}$. And from Lemma~\ref{hatfunctoriality}, $\mathcal{L}\mapsto\hat{\mathcal{L}}$ is an endofunctor on \dva-frames and skeletal morphisms. 
So it remains to check that this corestricts to double negation frames, and that the forgetful functor is adjoint.

Suppose $\mathcal{L}$ is corrigible. Then $a\in\hat{L}_-$ by definition means that $b\sqsubseteq_- a$ implies $b\leq a$. So
$a^{\bullet\bullet}\leq a$ in $\mathcal{L}$. But in $\hat{\mathcal{L}}$, $\con$ is just the restriction of $\con_{\mathcal{L}}$. 
So $a\mapsto a^\bullet$ in $\mathcal{L}$ restricts and corestricts to $\hat{\mathcal{{L}}}$.

The fact that the inclusion functor is adjoint to $\mathcal{L}\mapsto\hat{\mathcal{L}}$ follows immediately from Lemma~\ref{doubleneglemma}.
\end{proof}

This coreflection may be extended to all  \dva-frames with skeletal morphisms by considering a subcategory that includes double negation \dva-frames, the objects of which satisfy a property that is essentially the \dva-frame analogue of a known property of distributive lattices.

Say that a \dva-frame is \emph{dually subfit} if the following two conditions hold:
\[a\nleq_- b \Rightarrow \exists c\in L_-,\varphi\in L_+, \varphi\mathrel{\cancel{\con}} c\wedge a \text{ and } \varphi \con c \wedge b\]
and
\[\varphi\nleq_+ \psi \Rightarrow \exists \theta\in L_+,c\in L_-, \varphi\wedge \theta\mathrel{\cancel{\con}} c \text{ and } \psi\wedge \theta \con c.\]
Consider a \dva-frame $(L,L,\con,\tot)$ where $\con$ is defined by $a\con b$ if and only if $a\wedge b = 0$.  
Then dual subfitness amounts to the requirement that $a\nleq b$ implies there exists from $d$ so that $a\wedge d \neq 0$ and $b\wedge d= 0$.
This property (of distributive lattices) has had many names. 
Wallman \cite{Wallman38} refers to this as \emph{the disjunctive property} essentially because he formulated the implication as a disjunction. 
In set theoretic forcing, a generalization of this property to partially ordered sets is referred to as \emph{separativity}.
The name, \emph{subfitness}, for the dual property in point-free topology is due to Isbell \cite{I72}. 

\begin{lemma}\label{dualsubfitlemma1}
A \dva-frame $\mathcal{L}$ is dually subfit if and only if $a\sqsubseteq_- b$ 
implies $a\leq_- b$ and $\varphi\sqsubseteq_+\psi$ implies $\varphi\leq_+ \psi$.
\end{lemma}
\begin{proof}
Clearly, if $\mathcal{L}$ is dually subfit, then $a\nleq_- b$ implies $a\nsqsubseteq_- b$, and likewise for $\leq_+$. 
Conversely, if $\mathcal{L}$ is not dually subfit, the witnesses of failure provide witnesses that $a\sqsubseteq_- b$ but $a\nleq_- b$ or $\varphi\sqsubseteq_+\psi$ but not $\varphi\nleq_+\psi$.
\end{proof}

\begin{lemma}\label{dualsubfitlemma2}
For any \dva-frame $\mathcal{L}$, the \dva-frame $\hat{\mathcal{L}}$ is dually subfit. Moreover, for 
any dually subfit frame $\mathcal{L}$ is isomorphic to $\hat{\mathcal{L}}$.
\end{lemma}

\begin{proof}
Recall from Lemma \ref{quotientsq} that
$\hat{\mathcal{L}}$ is obtained precisely by taking the quotients of $L_-$ and $L_+$ for which $\sqsubseteq_-$ in $L_-$ becomes $\leq_-$ in $\hat{L}_-$, 
and $\sqsubseteq_+$ in $L_+$ becomes $\leq_+$ in $\hat{L}_+$.
And since $\con_{\hat{\mathcal{L}}}$ in $\hat{\mathcal{L}}$ is the restriction of $\con_{\mathcal{L}}$, $\hat{\mathcal{L}}$ satisfies the condition in Lemma~\ref{dualsubfitlemma1}.

For a dually subfit \dva-frame $\mathcal{L}$, $\sqsubseteq_-$ and $\leq_-$ coincide. So taking the quotient of $L_-$ with respect to $\sqsubseteq_-$ is trivial. The same is true in $L_+$.
\end{proof}

\begin{theorem}\label{dualsubfitCoreflection}
  The category of dually subfit \dva-frames is co-reflective in the category of all \dva-frames with skeletal morphisms.
\end{theorem}

\begin{proof}
The preceding lemmas do most of the work.
We only need to remark that if $\mathcal{L}$ is dually subfit, then any \dva-frame morphism $f\colon\mathcal{L}\to\mathcal{M}$ is skeletal because this simply means that $\leq_-$ and $\leq_+$ are preserved by $f_-$ and $f_+$.
\end{proof}

Theorems~\ref{doublenegtionCoreflection} and \ref{dualsubfitCoreflection} point toward several interesting research lines.
In fact, all the results of this section pertain to properties of $\con$, with no mention of $\tot$. 
Observe that that axioms mentioning $\con$ only amount to saying that $\con$ is a dcpo with bottom element $(0,0)$ and a distributive lattice when regarded as a sublattice of $L_+\times L_-^\partial$ (the order dual of $L_-$).
In fact, a very efficient description of the orders that can arise as $\con_{\mathcal L}$ for some \dva-frame $\mathcal{L}$ is that they are precisely the bounded distributive lattices in the category of dcpos with least element.
The objects have two cooperating orders:
the \emph{information order} induced by $L_+\times L_-$ is a dcpo with least element; the \emph{logical order} induced by $L_+\times L_-^\partial$ is a bounded distributive lattice.

On this view, the results of this section pertain to a category in which logic and information interact.
Pseudocomplementation needs re-interpretation.
Lemma~\ref{propspseudo} shows that $\varphi\con a$ implies $a^\bullet\con \varphi^\bullet$.
So $(-)^\bullet$ can be taken as an operation on $\con$ itself, sending $(a,\varphi)$ to the set of elements informationally disjoint from $(a,\varphi)$.
This suggests interesting questions of how logic and information interact.
For example, the law of double negation takes on a new meaning since pseudocomplementation has to do with information.
One result of this re-interpretation is, for example, that the law of double negation is strictly weaker than the law of excluded middle (which naturally refers to $\tot$).
Dual subfitness is still weaker than the law of double negation, but as of this writing we do not know how to interpret dual subfitness in terms of logic and information.
We expect that investigations into the \emph{logico-informational} interpretation of \dva-frames will bear fruit.

\section{Conclusion}

This paper addresses dense sub-objects for point-free bitopology in terms of \dva-frames. We have characterized extremal epimorphisms in general, and shown that a smallest dense one always exists.
Interestingly, the latter part of Isbell's theorem fails, in so far as the resulting \dva-frames are not typically Boolean. 
On the other hand, they do satisfy a weak law of separation (dual subfitness), and when pseudo-complements are well-behaved, a law of double negation. 

\Dva-frames provide only one formulation of point-free bitopology. Biframes offer an alternative, and are far better represented in the literature. 
Nevertheless for biframes, extremal epimorphisms and the analogue of Isbell's theorem have not been yet conclusively investigated (see \cite{FS11,PP14}). In the sequel to this paper, we turn our attention to these topics for biframes.

\section*{Acknowledgements}
Ale\v s Pultr's unflagging enthusiasm for mathematics and affection for fellow mathematicians is a model we find pleasure in emulating. 
His capacity to ask the compelling question is often quite uncanny. But of course, it comes from a depth of insight into where those questions may lead.

There is an old joke among modal logicians: "I wish to thank my host, without whom this would not be possible. I also wish to thank my family, without whom this would not be necessary.” We thank Ale\v s, without whom our mathematical lives would not have been possible. We also wish to thank him for his questions, without which this paper would not have been necessary.

The authors acknowledge financial support from the grant MTM2015-63608-P (MINECO/FEDER, UE) and from the Basque Government (grant IT974-16 and Postdoctoral Fellowship POS\_2016\_2\_0032).

\end{document}